\documentclass[11pt]{article}
\usepackage{amsmath}
\usepackage{mathrsfs}
\usepackage{amsfonts}

\usepackage{amsfonts, amsmath, amssymb}
\usepackage{amssymb,amsfonts,amsmath,
latexsym, epsfig,cite, psfrag,eepic,colordvi}
\usepackage{amscd,graphics}
\usepackage{lineno}

\newcommand{\qed}{\hfill $\Box $}

\newtheorem{theorem}{Theorem}[section]

\newtheorem{coro}[theorem]{Corollary}

\begin{document}

\title{Simplified existence theorems on all fractional $[a,b]$-factors
\thanks{This work was supported
 by  the National Natural
Science Foundation of China No. 11101329 and the Fundamental
Research Funds for the Central Universities.} }

\author{Hongliang Lu\thanks{Corresponding email:
luhongliang215@sina.com (H. Lu)}
\\ {\small Department of Mathematics}
\\ {\small Xi'an Jiaotong University, Xi'an 710049, PR China}
}

\date{}

\maketitle

\date{}

\maketitle

\begin{abstract}

Let $G$ be a graph with order $n$ and let $g, f : V (G)\rightarrow
 N$ such that $g(v)\leq f(v)$ for all $v\in V(G)$. We say that $G$
has all fractional $(g, f )$-factors if $G$ has a fractional
$p$-factor for every $p: V (G)\rightarrow N$ such that $g(v)\leq
p(v)\leq f (v)$ for every $v\in V(G)$. Let $a<b$ be two positive
integers. 
If $g\equiv a$, $f\equiv b$ and $G$ has all fractional
$(g,f)$-factors, then   we say that  $G$ has all fractional
$[a,b]$-factors.  Suppose that $n$ is sufficiently large  for $a$
and $b$.

This paper contains two results on the existence of all
$(g,f)$-factors of  graphs. First, we derive from Anstee's
fractional $(g,f)$-factor theorem a similar characterization
 for the property
of all fractional $(g,f)$-factors. Second, we show that $G$ has all
fractional $[a, b]$-factors  if the minimum degree is at least
$\frac{1}{4a}((a+b-1)^2+4b)$ and every pair of nonadjacent vertices
has cardinality of the neighborhood union at least $bn/(a + b)$.
These lower bounds are sharp.

\end{abstract}


\section{Introduction}

We consider finite undirected graphs without loops or multiple
edges. Let $G$ be a graph with vertex set $V (G)$ and edge set
$E(G)$. For $x \in V (G)$, we denote by $d_G(x)$ the degree of $x$
in $G$ and by $N_G(x)$ the set of vertices adjacent to $x$ in $G$.
We write $N_G[x]$ for $N_G(x)\cup\{x\}$. The minimum degree and the
maximum degree of $G$ are denoted by $\delta(G)$ and $\Delta(G)$,
respectively. For a subset $S\subseteq V (G)$, let $N_G(S)$ denote
the union of $N_G(x)$ for every $x\in S$ and $G[S]$ is the subgraph
of $G$ induced by $S$. We write $G-S$ for $G[V (G)\setminus S]$ and
$f(U)$ for $\sum_{v\in U} f (u)$.  For $M\subseteq E(G)$, let $G[M]$
denote the subgraph of $G$ induced by $M$.
 For two disjoint vertex subsets
$A$ and $B$ of $G$, the number of edges joining $A$ to $B$ is
denoted by $e_G(A,B)$. For a vertex $x\in V(G)$ and an edge $e\in
E(G)$,   we write $x\sim e$ if $x$ is incident with $e$. The
\emph{join} $G + H$ denotes the graph with vertex set $V (G)\cup V
(H)$ and edge set
$$E(G + H) = E(G)\cup E(H)\cup \{xy\ |\  x \in V (G)\ \mbox{and}\ y\in V(H)\}.$$

Let $g$ and $f$ be two integer-valued functions defined on $V(G)$
such that $0\leq g(x)\leq f(x)$ for all $x\in V(G)$. A
\emph{$(g,f)$-factor} of $G$  is a spanning subgraph $F$ of $G$
satisfying $g(x)\leq d_F(x)\leq f(x)$ for all $x\in V(G)$. Let $a<b$
be two integers. A $(g, f )$-factor is called an \emph{$[a,
b]$-factor} if $g(x)\equiv a$ and $f (x)\equiv b$. Let $h :
E(G)\rightarrow [0, 1]$ be a function. If $g(x)\leq \sum_{x\sim
e}h(e)\leq f(x)$ holds for any vertex $x \in V (G)$, then we call
graph $F$ with vertex set $V(G)$ and edge set $E_h$ a
\emph{fractional $(g,f)$-factor} of $G$ with indicator function $h$,
where $E_h =\{e \in  E(G)\ |\ h(e)>0\}$.  If $f(v)=g(v)$ for all
$v\in V(G)$, then a fractional $(g,f)$-factor is called
\emph{fractional $f$-factor}.
If $h(e)\in \{0,1\}$ for every edge $e\in E(G)$, then $F$ is just a
$(g,f)$-factor of $G$. Clearly, if $G$ contains a $(g,f)$-factor,
then it also contains a fractional $(g,f)$-factor.


Lov\'asz \cite{Lov70} gave a characterization of graphs having a
$(g, f )$-factor.  More generally, one could consider the existence
of $p$-factors, where $p: V(G)\rightarrow Z^+$ is a function such
that $g(v)\leq p(v)\leq f (v)$ for every $v\in V(G)$ and $p(V)\equiv
0$ (mod 2). We say that $G$ has all \emph{$(g, f )$-factors}
  if  $G$ has a
$p$-factor for every $p$ described above.   Kano and Tokushige
\cite{Kano92} gave a sufficient condition for a graph to have all
$(g,f)$-factors in terms of the minimum degree.  Niessen
\cite{Niessen} obtained a necessary and sufficient condition for a
graph to have all $(g,f)$-factors.

\begin{theorem}[Kano and Tokushige, \cite{Kano92}]\label{Kano}
Let $G$ be a connected graph of order $n$, $a$ and $b$ be integers
such that $1\leq a\leq b$ and $f: V(G)\rightarrow
\{a,a+1,\ldots,b\}$. Suppose that $n\geq (a+b)^2/a$ and
$f(V(G))\equiv 0 \pmod 2$. If $\delta(G)>\frac{an-2}{a+b}$, then $G$
contains an $f$-factor.

\end{theorem}

\begin{theorem}[Niessen, \cite{Niessen}] A graph
$G$ has all $(g, f )$-factors if and only if
\[ g(S)-f(T)+\sum_{x\in T}d_{G-S}(x)-q^*_G(S, T, g, f )\geq\left\{\begin{array}{ll}
-1, &\text{if }\
 f\neq g, \\[5pt]
0, &\text{if } f=g.
\end{array}\right.
\]
for all disjoint sets $S, T\subseteq V$, where $q^*_G(S, T, g, f )$
denotes the number of components $C$ of $G-(S\cup T)$ such that
there exists a vertex $v \in V(C)$ with $g(v)< f (v)$ or $e_G(V(C),
T)+ f (V(C))\equiv 1 \pmod 2$.

\end{theorem}

   Anstee \cite{Anstee}
gave the following necessary and sufficient condition for a graph to
have a fractional $(g,f)$-factor.
\begin{theorem}[Anstee, \cite{Anstee}]\label{Anstee}
Let $G$ be a graph and $g,f: V(G)\rightarrow Z^+$ be two integer
functions such that $g(v)\leq f(v)$ for all $v\in V(G)$. Then $G$
has a fractional $(g, f)$-factor if and only if for any subset
$S\subseteq V(G)$, we have
\begin{align*}
f(S)-g(T)+ \sum_{v\in T} d_{G-S}(v)\geq 0,
\end{align*}
where $T = \{v\ |\  v \in V (G)-S \mbox{ and } d_{G-S}(v)<g(v)\}$.
\end{theorem}

More naturally, we consider a similar property of graphs having all
$(g,f)$-factors.
 Let $q$ be an
integer-valued function defined on $G$ such that $g(v)\leq q(v)\leq
f (v)$ for every $v\in V(G)$. We  say that $G$ has all
\emph{fractional $(g, f )$-factors}
  if  $G$ has a fractional
$q$-factor for every $q$ described above. If $g\equiv a$, $f\equiv
b$ and $G$ has all fractional $(g,f)$-factors, then   we say that
$G$ has \emph{all fractional $[a,b]$-factors}.

The following result is obtained by Zhang \cite{Zhang}.
\begin{theorem}[Zhang, \cite{Zhang}]\label{Zhang}
If $G$ has a fractional $(g, f)$-factor, then it must have a
fractional $(g, f)$-factor $F$ with indicator function $h$  such
that $h(e)\in \{0, 1/2 , 1\}$ for every edge $e \in E(G)$.

\end{theorem}

 Let $G^*$ be  a graph obtained from $G$ by replacing every edge by
two multiple edges. By Theorem \ref{Zhang}, $G$ has a fractional
$f$-factor if and only if $G^*$ has a $(2f)$-factor.  Hence $G$ has
all fractional $(g,f)$-factors if and only if $G^*$ has a
$(2p)$-factor, for every integer-valued function $p$ defined $G$
which satisfies   $g(v)\leq p(v)\leq f(v)$ for all $v\in V(G)$.

In this paper, we first  give
a  necessary and sufficient condition for a graph to have all
fractional $(g,f)$-factors. 

\begin{theorem}\label{lu}
Let $G$ be a graph and $g,f: V(G)\rightarrow Z^+$  be two integer
functions such that $g(x)\leq f(x)$ for all $x\in V(G)$. Then $G$
has all fractional $(g, f)$-factors if and only if for any subset
$S\subseteq  V(G)$, we have
$$g(S)-f(T)+\sum_{x\in T}d_{G-S}(x)\geq 0
$$ where $T = \{v\ |\  v \in V (G)-S \mbox{ and } d_{G-S}(v)<f(v)\}$.
\end{theorem}

If $g\equiv a$ and $f\equiv b$, then  by Theorem \ref{lu},  we
obtain the following result.
\begin{coro}\label{lu2}
Let $G$ be a graph and $a<b$ be two positive integers.
 Then $G$ has all fractional $[a, b]$-factors  if and only if
for any subset $S\subseteq  V(G)$, we have
$$a|S|-b|T|+\sum_{x\in T}d_{G-S}(x)\geq 0,
$$ where $T = \{v\ |\  v \in V (G)-S \mbox{ and } d_{G-S}(v)<b\}$.
\end{coro}
%
%
%
By  Corollary \ref{lu2}, we obtain a sufficient condition for a
graph to have all  fractional $[a,b]$-factors, which is closely
related with Theorem \ref{Kano}.

\begin{theorem}\label{lu3}
Let $a<b$ be two positive integers. Let $G$ be a graph with order
$n\geq 2(a+b)(a+b-1)/a$ and minimum degree $\delta(G)\geq
\frac{1}{4a}((a+b-1)^2+4b)$. If $|N_G(u)\cup N_G(v)| \geq
\frac{bn}{a+b} $ for any two nonadjacent vertices $u$ and $v$ in
$G$, then $G$ has all fractional $[a, b]$-factors.

\end{theorem}

%
%
%

\section{The proof of Theorems \ref{lu} and \ref{lu3}}

\noindent\textbf{Proof of Theorem \ref{lu}.} We first prove the
sufficiency. Let $h: V(G)\rightarrow Z^+$ be an arbitrary integer
function such that $g(v)\leq h(v)\leq f (v)$ for every $v\in V(G)$.
For any $S\subseteq V(G)$, let $T = \{v\ |\ v \in V (G)-S \mbox{ and
} d_{G-S}(v)<f(v)\}$ and $T'=\{x\ |\ x \in V (G)-S \mbox{ and }
d_{G-S}(x)<h(x)\}$, we have
\begin{align*}
h(S)-h(T')+\sum_{v\in T'}d_{G-S}(v)\geq g(S)-f(T)+\sum_{v\in
T}d_{G-S}(v) \geq 0.
\end{align*}
By Theorem \ref{Anstee}, $G$ has a fractional $h$-factor.

Now we prove the necessary. Conversely, suppose that there exists
$S\subseteq V(G)$ such that $$g(S)-f(T)+\sum_{x\in T}d_{G-S}(x)<0,
$$
where $T = \{v\ |\ v \in V (G)-S \mbox{ and } d_{G-S}(v)<f(v)\}$.
Let $h(v)=g(v)$ for all $v\in S$, $h(u)=f(u)$ for all $u\in V(G)-S$.
Then we have
\begin{align*}
0&>g(S)-f(T)+\sum_{x\in T}d_{G-S}(x)\\
&=h(S)-h(T)+\sum_{x\in T}d_{G-S}(x).
\end{align*}
So,   $G$ contains no fractional $h$-factors by Theorem
\ref{Anstee}, a contradiction. This completes the proof. \qed

%
%
%

\noindent\textbf{Proof of Theorem \ref{lu3}.} Suppose that the
result doesn't hold. By Corollary \ref{lu2}, there exists
$S\subseteq V(G)$ such that
\begin{align}\label{eq:1}
a|S|-b|T|+\sum_{x\in T}d_{G-S}(x)< 0,
\end{align}
where $T = \{v\ |\  v \in V (G)-S \mbox{ and } d_{G-S}(v)<b\}$. Let
$s=|S|$ and $t=|T|$. Note that $|T| \neq \emptyset$, otherwise,
$$a|S|-b|\emptyset|+\sum_{x\in \emptyset}d_{G-S}(x)=a|S|\geq 0,$$
contradicting to (\ref{eq:1}). 
For a vertex $x\in T$,  we write $N_T(x)=N_G(x)\cap T$ and
$N_T[x]=(N_G(x)\cap T)\cup \{x\}$.

Now
 we define
\begin{align*}
h_{1}=\min \{d_{G-S}(v)\ |\ v\in T\}.
\end{align*}
And let $x_1$ be a vertex in
$T$ satisfying $d_{G-S}(x_1) = h_1$. Further, if $T-N_{T}[x_1]\neq
\emptyset$, we define
\begin{align*}
h_{2}=\min \{d_{G-S}(v)\ |\ v\in T-N_{T}[x_1]\},
\end{align*}
and let $x_2$ be a vertex in $T - N_{T}[x_1]$ satisfying
$d_{G-S}(x_2) = h_2$. Then we have $h_1\leq  h_2\leq b- 1$ and
$|N_G(x_1)\cup N_G(x_2)|\leq s+h_1+h_2$.


Now in order to prove the correctness of the theorem, we will deduce
some contradictions according to the following two cases.

\medskip
\textbf{ Case 1.~} $T = N_{T}[x_1]$.
\medskip

Since $\delta(G)\geq \frac{1}{4a}((a+b-1)^2+4b)$, $t\leq h_1+1\leq
b$ and $h_1+s\geq d_G(x_1)\geq \delta(G)$, by (\ref{eq:1}) and the
definition of $h_1$, we have
\begin{align*}
0 &\geq as-bt+\sum_{x\in T}d_{G-S}(x)+1\\
&\geq as-bt + h_{1}t + 1\\
&\geq a(\delta(G)-h_1)-(b-h_{1})t+1\\
&\geq a(\delta(G)-h_1)-(b-h_{1})(h_1+1)+1\\
&=h_1^2-(a+b-1)h_1+a\delta(G)-b+1\\
&\geq (h_1-\frac{a+b-1}{2})^2+1\geq 1,
\end{align*}
a contradiction.

\medskip
\textbf{ Case 2.~} $T - N_{T}[x_1]\neq \emptyset$.
\medskip

Let $W=G-S-T$ and $w$ be the number of components of $W$. Then
$w\leq n-s-t$. Let $p=|N_{T}[x_1]|$, then $h_1\geq p-1$. And by $T -
N_{T}[x_1]\neq \emptyset$, we have $t\geq p+1$. Since $0\leq w\leq
n-s-t$, $h_1-h_2\leq 0$, $b-h_{2}\geq 1$, and (\ref{eq:1}), we have
\begin{align*}
(n-s-t)(b-h_{2})&\geq as-bt+\sum_{x\in T}d_{G-S}(x)+1\\
&\geq as-bt+h_{1}p+h_{2}(t-p)+1\\
&=as+(h_{1}-h_2)p+(h_{2}-b)t+1\\
&\geq as+(h_{1}-h_2)(h_1+1)+(h_{2}-b)t+1.
\end{align*}
Hence,
\begin{align*}
(n-s)(b-h_{2})\geq as+(h_1-h_2)(h_1+1)+1.
\end{align*}
Then it follows from the inequality above that
\begin{align}\label{eq:2}
n(b-h_{2})- (a+b-h_2)s+(h_2-h_1)(h_1+1)-1\geq 0.
\end{align}

 Since two vertices $x_1$ and $x_2$ are not adjacent, by the
condition of the theorem, the following inequalities hold:
\begin{align*}
\frac{bn}{a+b}\leq |N_{G}(x_1)\cup N_G(x_2)|\leq h_1+h_2+s,
\end{align*}
which implies
\begin{align}\label{eq:6}
s\geq \frac{bn}{a+b}-(h_1+h_2).
\end{align}
Using (\ref{eq:2}) and (\ref{eq:6}), $h_1\leq h_2\leq b-1$, and
$n\geq 2(a + b)(a + b - 1)/a$, we obtain
\begin{align*}
0&\leq n(b-h_{2})-(a+b-h_2)s+(h_2-h_1)(h_1+1)-1\\
&\leq
n(b-h_{2})-(a+b-h_2)(\frac{bn}{a+b}-(h_1+h_2))+(h_2-h_1)(h_1+1)-1\\
&=\frac{-anh_2}{a+b}+(a+b+1)h_2-h_2^2+(a+b-1)h_1-h_1^2-1\\
&\leq \frac{-anh_1}{a+b}+(a+b+1)h_1-h_1^2+(a+b-1)h_1-h_1^2-1\\
&=\frac{-anh_1}{a+b}+2(a+b)h_1-2h_1^2-1\\
&\leq 2h_1-2h_1^2-1\leq -1,
\end{align*}
also a contradiction. This last contradiction completes our proof.
\qed


%

\noindent\textbf{Remark:}  In Theorem \ref{lu3}, the bound in the
assumption
\begin{align*}
|N_G(u)\cup N_G(v)|\geq  \frac{bn}{a+b}
\end{align*}
is best possible. We can show this by constructing a  graph $G=
K_{bm}+(am+1)K_{1}$, where $m$ is any positive integer. In this
graph $G$, there exist two nonadjacent vertices $u$ and $v$ which
degrees are equal to the minimum degree of value $bm$. Obviously, we
know
\begin{align*}
\frac{bn}{a+b}=\frac{b(am+bm+1)}{a+b}>bm=|N_G(u)\cup
N_G(v)|>\frac{bn}{a+b}-1.
\end{align*}
Let $h:V(G)\rightarrow Z^+$ such that   $h(u)=a$ for all  $u\in
V(K_{bm})$ and $h(v)=b$ for all $v\in V(G)-V(K_{bm})$.  Then we have
\begin{align*}
a|S|-b|T|+\sum_{x\in T}d_{G-S}(x)=a|S|-b|T|<0,
\end{align*}
where  $S=V(K_{bm})$ and $T=V(G)-V(K_{bm})$. So by Theorem \ref{lu},
 graph $G$ contains no fractional $h$-factors.

In the following, we show that the bound on minimum degree in
Theorem \ref{lu3} is also best possible. Let $1\leq a<b$ be two
integers such that $a+b$ is odd and $r$ be a sufficient large
integer. Let $m$ be the maximum integer such that
$m<\frac{1}{4a}((a+b-1)^2+4b)-\frac{(a+b-1)}{2}$. Let
$G=K_m+(K_r\cup K_{(a+b+1)/2})$.
 Then we have
$\delta(G)=m+\frac{(a+b-1)}{2}<\frac{1}{4a}((a+b-1)^2+4b)$ and
$|N_{G}(x)\cup N_G(y)|\geq \frac{bn}{a+b}$ for any two nonadjacent
vertices. Let $h:V(G)\rightarrow Z^+$ be an integer-valued function
such that $h(u)=a$ for  all $u\in V(K_{m})$ and $h(v)=b$ for all
$v\in V(G)-V(K_{m})$. Similarly,  by Theorem \ref{lu2},  it is easy
to see that this graph $G$ contains  no fractional $h$-factors.


%

\end{document}